\documentclass[11pt]{article}

\usepackage{amsthm,amssymb,amscd}
\theoremstyle{plain}
\newtheorem{rem}{Remark}[section]

\newtheorem{cor}{Corollary}[section]

\newtheorem{thm}{Theorem}[section]
\theoremstyle{definition}

\newcommand{\bb}[1]{\mbox{$\mathbb{#1}$}}

\title{A generalized Verdier-type Riemann-Roch theorem for Chern-Schwartz-MacPherson
classes}
\author{J\"{o}rg Sch\"{u}rmann\thanks{Westf. Wilhelms-Universit\"{a}t, SFB 478
"Geometrische Strukturen in der Mathematik",
Hittorfstr.27, 48149 M\"{u}nster, Germany,
E-mail: jschuerm@math.uni-muenster.de}
\thanks{Partially supported by the EU Research Training Network ``Geometric Analysis''}}
\date{ }

\begin{document}
\bibliographystyle{plain}

\maketitle

\begin{abstract}
In this paper we give a general formula for the defect appearing in the Verdier- 
type Riemann-Roch formula for Chern-Schwartz-MacPherson classes in the case of a regular embedding
(and for suitable local complete intersection morphisms).
Our proof of this formula uses the "constructible function version" of Verdier's specialization
functor $Sp_{X\backslash Y}$ (for constructible (complexes of) sheaves), together with  a specialization
 property of 
Chern-Schwartz-MacPherson classes and the corresponding Riemann-Roch theorem for smooth
morphisms. As a very special case we get a general formula for the "Milnor-class" of a local
complete intersection in a smooth manifold (generalizing many recent results of different authors).
Our formula becomes very simple in the case of a regular embedding of codimension one,
and gives a far reaching generalization of Parusi\'{n}ski-Pragacz's formula for the
"Milnor-class" of a hypersurface in a (compact) smooth manifold. 
\end{abstract}

\section*{Introduction}  
The Chern-Schwartz-MacPherson class $c_{*}$ is the unique transformation from constructible
functions to homology, commuting with proper push forward and satisfying for smooth $X$
the normalization
\begin{equation}
c_{*}(1_{X}) = c^{*}(TX) \cap [X] \;,
\end{equation}
with $c^{*}(TX)\cap$ the action of the total Chern-class of the tangent-bundle $TX$ on the
fundamental-class $[X]$ of $X$. Note that the uniqueness follows from resolution of singularities
and the existence of such a transformation was conjectured by Deligne and Grothendieck
(compare \cite[p.361-368]{Gr}).
It was first constructed by MacPherson in the complex algebraic context in his classical paper
\cite{M}, using his celebrated "Euler-obstruction" and "graph-construction".
In his origional work MacPherson used as a suitable homology theory the Borel-Moore homology
$H^{BM}_{*}(\cdot,\bb{Z})$ of the corresponding analytic space. But his proof also applies to
the case, where one uses the algebraic Chow-homology theory $A_{*}$ of Fulton \cite{Ful}
(compare \cite[Ex. 19.1.7]{Ful}). Moreover both approaches are linked by the cycle map
$cl : A_{*} \to H^{BM}_{*}(\cdot,\bb{Z})$ \cite[Chapter 19]{Ful}.\\

The analyticity of the graph construction was finally solved by Kwieci\'{n}ski in his
thesis \cite{Kw} so that MacPherson's proof also applies to the complex analytic context.
Moreover, for complex spaces which can be embedded into complex manifolds, the
MacPherson Chern-class $c_{*}(1_{X})$ of $X$ is isomorphic via Alexander duality 
(by work of Brasselet-Schwartz \cite[Exp. 6]{CEP}) to the classes introduced already before
by Schwartz \cite{Schw}. Finally a new approach to these Chern-classes (in the embeddable context)
comes from the "micro-local viewpoint" through the theory of Lagrangian cycles
\cite{Fu2,Gi,Sab,Sch2}. This theory was originally described in terms of characteristic
cycles of holonomic D-modules or constructible (complexes of) sheaves, but only depends on the
associated constructible functions (For a direct approach to Lagrangian cycles in terms of
constructible functions compare with \cite{Fu2,Sch3}). Moreover, the Lagrangian approach
has been extended in the algebraic context to (embeddable) schemes of finite type over
a ground field of characteristic zero \cite{Ken}.\\

The Verdier-type Riemann-Roch formula one is looking for is the commutativity of the following 
diagram:
\begin{equation} \label{Verdier} \begin{CD}
F(Y) @> c_{*} >> H_{*}(Y)\\
@V f^{*} VV  @VV c^{*}(T_{f})\cap f^{*} V\\
F(X) @>> c_{*} > H_{*}(X) 
\end{CD} \end{equation}
for a (local complete intersection) morphism $f: X\to Y$, which has a global factorization
into a regular embedding $i: X\to Z$ followed by a smooth morphism $p: Z\to Y$.
Here we use the following notations:
\begin{itemize}
\item $f^{*}: F(Y)\to F(X);\; \alpha\mapsto \alpha \circ f$ is the usual pullback of constructible functions.
\item $H_{*}(\cdot)$ the relevant homology theory $A_{*}(\cdot)$ (or $H^{BM}_{*}(\cdot,\bb{Z})$)
depending on the algebraic (or complex analytic) context.
\item $c^{*}(T_{f})\cap := \Bigl( c^{*}(i^{*}T_{p})\cup \bigl(c^{*}(N_{X}Z)\bigr)^{-1}\Bigr) 
\cap \;,$ with $T_{p}$ the relative tangent bundle of the smooth morphism $p$
and $N_{X}Z$ the normal bundle of the regular embedding $X\to Z$.
\item Finally $f^{*}: H_{*}(Y)\to H_{*}(X)$ is the Gysin transformation for such a
morphism, i.e. the composition of the smooth pullback $p^{*}$ and the Gysin transformation
$i^{*}$ for the regular embbeding $i$ (which is defined with the help of the deformation to the
normal cone $N_{X}Z$). For the construction of $i^{*}$ in the algebraic (or analytic) context 
compare with \cite[Ch.5,Ch.6]{Ful} (or \cite[Exp.IX]{GA}).
\end{itemize}

Note that $c^{*}(T_{f})\cap$ and $f^{*}: H_{*}(Y)\to H_{*}(X)$ are independent of the chosen
factorization $f=p\circ i$ of $f$.
This commutativity is motivated by a corresponding result of Verdier \cite[p.214, thm.7.1]{GA}
for the "Todd-homology transformation" $\tau_{*}$ of Baum-Fulton-MacPherson
(compare with \cite[thm.18.2(3)]{Ful} and especially also with \cite[p.11,p.122]{FM} 
for the origin of the name for such a formula). \\

A similar Verdier-Riemann-Roch theorem 
for Chern-Schwartz-MacPherson classes $c_{*}$ is true
for a smooth morphism $p$ (\cite[thm.2.2]{Y} and compare also with \cite[p.111]{FM}), i.e. the following
diagram is commutative for a smooth morphism $p: X\to Y$:

\begin{equation} \label{smoothVerdier} \begin{CD}
F(Y) @> c_{*} >> H_{*}(Y)\\
@V p^{*} VV  @VV c^{*}(T_{p})\cap p^{*} V\\
F(X) @>> c_{*} > H_{*}(X)\;, 
\end{CD} \end{equation}
with $T_{p}$ the relative tangent bundle of $p$.
Note that Yokura works in the complex algebraic 
context, but his simple proof applies also to the analytic (or algebraic context over a ground field
of characteristic zero), since it uses resolution of singularities, a smooth base-change
property (for constructible functions and for the corresponding homology theory), together with the
defining properties (normalization, proper push forward) of the Chern-classes $c_{*}$
(and for a Lagrangian approach compare with \cite{Sch2}).\\

But in general the diagram (\ref{Verdier}) is not commutative!
The factorization
\begin{equation} \label{factVerdier} \begin{CD}
F(Y) @> c_{*} >> H_{*}(Y)\\
@V p^{*} VV  @VV c^{*}(T_{p})\cap p^{*} V\\
F(Z) @> c_{*} >> H_{*}(Z)\\
@V i^{*} VV  @VV c^{*}(N_{X}Z)^{-1}\cap i^{*} V\\
F(X) @>> c_{*} > H_{*}(X) 
\end{CD} \end{equation}

implies that the problems come from the case of a regular embedding $i: X\hookrightarrow Y$.
For the special case of a smooth ambient space $Y$
\[\bigl(c^{*}(N_{X}Y)\bigr)^{-1} \cap i^{*}(c_{*}(1_{Y}))
= \bigl(c^{*}(N_{X}Y)\bigr)^{-1} \cap c^{*}(i^{*}TY) \cap [X]\] 
is the Fulton-Johnson Chern class $c^{FJ}(X)$ of $X$ \cite[Ex.4.2.6]{Ful}, which is in
general different from the Chern-class $c_{*}(1_{X}) = c_{*}(i^{*}1_{Y})$ of $X$
(compare \cite[Ex.3.1]{Y}). The difference class 
\begin{equation} \label{Milnor}
M_{0}(X) := c^{FJ}(X)
-  c_{*}(1_{X})
\end{equation}
is called the "Milnor class" of $X$ (up to suitable sign-conventions) and has been recently 
studied in some cases by many authors using quite different methods.
Compare for example with \cite{Al,BLSS,BLSS2,PP,Su,OY,Y3} and especially with the nice survey 
articles \cite{Br,Y2} of Brasselet and Yokura. We will explain the relations to our results in the end of 
this work.
Note that these authors use a different sign-convention for their Milnor-class $M(X)$,
and for this reason we use the notation $M_{0}(X)$ (following \cite{Y3}).
Moreover, in most cases they consider (only) the case of a "global" complete intersection $X$ in $Y$,
i.e. $X$ is defined by a cartesian diagram
\[ \begin{CD}
X @>i >>  Y \\
@VVV @VV s V \\
Y @>> 0  > E\;,
\end{CD} \]
with $s$ a section of a vector-bundle $E\to Y$ which is generically transverse to the zero-section
$0: Y\to E$. So in this case one gets $N_{X}Y = E|X$.\\

In this paper we explain in a very general context the non-commutativity of (\ref{Verdier})
in the case of a regular embedding  $i: X\hookrightarrow Y$ 
using the "constructible function version" of Verdier's specialization
to the normal cone (i.e. of Verdier's specialization transformation $Sp_{X\backslash Y}$ 
\cite{V} for constructible (complexes of) sheaves):
\begin{equation}\label{special}
Sp_{X\backslash Y}: F(Y) \to F_{mon}(C_{X}Y)\;,
\end{equation}
with $F_{mon}(C_{X}Y)$ the group of "monodromic" constructible functions 
on the normal cone (i.e. constructible 
functions on $C_{X}Y$, which are invariant under the natural $\bb{C}^{*}$-action on this cone
over $X$). Note that this definition works for any closed embedding 
$i: X\hookrightarrow Y$ of a subspace (with $C_{X}Y = N_{X}Y$ the normal bundle in the case of a 
regular embedding).
We explain in the next chapter the definition and basic
properties of $Sp_{X\backslash Y}$ (following \cite{V}). Let us only mention 
that $Sp_{X\backslash Y}(\alpha)|X = \alpha|X$ for any $\alpha\in F(Y)$, where we consider $X$
 as a subspace
of the normal cone $C_{X}Y$ as the "zero- or vertex-"section $k: X\hookrightarrow C_{X}Y$.
Let $\pi: C_{X}Y\to X$ be the projection. \\

Then we introduce the "generalized vanishing cycles"
$\Phi_{i}$ of $i$ as the transformation
\begin{equation}\label{vanish}
\Phi_{i}:= Sp_{X\backslash Y}
- \pi^{*}i^{*} 
: F(Y) \to F_{mon}(C_{X}Y)\;.
\end{equation}
Note that $\Phi_{i}(\alpha)$  is a monodromic constructible function, which is zero on $X$
(for any $\alpha\in F(Y)$). The notation "generalized vanishing cycles" is motivated by the 
following observation (compare \cite[(Sp6) Normalisation,p.353]{V} and the next chapter):
Assume $i: X\hookrightarrow Y$ is a regular embedding of codimension one (i.e. a hypersurface on
$Y$ locally given by one equation $\{g=0\}$). By the above properties of $\Phi_{i}$
(and the fact that $\pi: C_{X}Y = N_{X}Y \to X$ is a line-bundle), there exists a welldefined
transformation $\mu: F(Y)\to F(X)$, with
\begin{equation} \label{defmu}
\Phi_{i} = (\pi^{*} - k_{*})\circ \mu \;.
\end{equation}
Here $k_{*}:  F(X)\to F(N_{X}Y)$ is just  extension by zero.
If in addition $X$ is (locally) given by the equation $\{g=0\}$ on $Y$, then
$\mu$ is just the "constructible function version" of  Deligne's vanishing cycle functor
$\phi_{g}$  for constructible (complexes of) sheaves.
Especially, this vanishing (and also the corresponding nearby) cycle functor is independent of the chosen
local equation of $X$.  This was one of the motivations of Verdier for introducing his specialization
functor $Sp_{X\backslash Y}$ (compare \cite[p.332]{V}).\\

Now we are ready to formulate the main result of this paper. For simplicity we work for the rest of this 
introduction in 
\begin{itemize}
\item the complex algebraic context with algebraic constructible functions $F(\cdot)$
on the associated complex analytic spaces of seperated schemes of finite type over 
$spec(\bb{C})$ and $H_{*}(\cdot)$ one of the two possible homology theories
(i.e. the Chow-homology $A_{*}(\cdot)$ or the Borel-Moore homology $H^{BM}_{*}(\cdot,\bb{Z})$
of the associated complex analytic space),
\item or in the complex analytic context with {\bf compact} analytic spaces.
\end{itemize} 
Then our main result is the following generalized Verdier-type Riemann-Roch theorem for
 Chern-Schwartz-MacPherson classes (asked for in \cite{Y,Y2,BY}): 
\begin{thm} \label{thm:main}
Consider a regular embedding  $i: X\hookrightarrow Y$
of complex "algebraic" (or compact analytic) spaces. Then one has for the "defect" in the
corresponding Verdier-type Riemann-Roch-formula the following explicit description:
\begin{equation} \label{defect}
c^{*}(N_{X}Y)^{-1} \cap i^{*}\bigl(c_{*}(\cdot)\bigr)
- c_{*}\bigl(i^{*}(\cdot)\bigr) 
= c^{*}(N_{X}Y)^{-1} \cap k^{*}\Bigl(c_{*}\bigl(\Phi_{i}(\cdot)\bigr)\Bigr)\;.
\end{equation}
Here $k: X\hookrightarrow N_{X}Y$ is the "zero-section" of the normal-bundle
$\pi: N_{X}Y\to X$ and $k^{*}$ on homology is the inverse 
of the isomorphism $\pi^{*}: H_{*}(X)\to H_{*}(N_{X}Y)$.
\end{thm}

The above "defect" formula follows by the definition of $\Phi_{i}$ directly from the following 
two results:
\begin{equation} \label{pullback}
c_{*}(\cdot) = c^{*}(N_{X}Y)^{-1} \cap k^{*}\Bigl(c_{*}\bigl(\pi^{*}(\cdot)\bigr)\Bigr)\;,
\end{equation}
which is a simple application of the corresponding Verdier Riemann-Roch theorem (\ref{smoothVerdier})
for the smooth
 projection $\pi: N_{X}Y\to X$:
\[ c_{*}\bigl(\pi^{*}(\cdot)\bigr) = c^{*}(T_{\pi}) \cap \pi^{*}\bigl( c_{*}(\cdot)\bigr) \;, \]
with $k^{*}\bigl(c^{*}(T_{\pi}) \cap \pi^{*}\bigr) = c^{*}(k ^{*}T_{\pi}) \cap$ (since $ k ^{*}T_{\pi}$
is isomorphic to $N_{X}Y$, and $k^{*}\pi^{*} = id$).
This part works without any compactness assumptions.\\

The other result is the formula
\begin{equation} \label{sp}
i^{*}\bigl(c_{*}(\cdot)\bigr) = k^{*}\Bigl(c_{*}\bigl(Sp_{X\backslash Y}(\cdot)\bigr)\Bigr) \;,
\end{equation}
which follows from the definition of the Gysin homomorphism $i^{*} = k^{*}\circ sp$ (based on the
deformation to the normal cone) and the specialization
property 
\begin{equation} \label{spec}
sp \circ c_{*} = c_{*} \circ  Sp_{X\backslash Y}
\end{equation}
of the Chern-class transformation $c_{*}$ (Compare with \cite[p.352]{Ful} and \cite[thm.6.5,p.211]{GA}
for the corresponding specialization property of
the "Todd-homology transformation" $\tau_{*}$ of Baum-Fulton-MacPherson).\\

This specialization property (\ref{spec}) of $c_{*}$ follows from a corresponding result
with respect to the "constructible function 
counterpart" of Deligne's nearby cycle functor (used in the definition of $ Sp_{X\backslash Y}$).
In the analytic context this goes back to Verdier  \cite[Exp.7, thm.5.1]{CEP}, and it is exactly 
at this point where we use in the analytic context our compactness assumption.
Other proofs of this specialization property for $c_{*}$ (in the embedded context) have
been given by the theory of Lagrangian cycles in \cite{Fu,Gi,Sab,Sch3}.
Moreover, this Lagrangian approach was extended by Kennedy \cite[thm.2)]{Ken2} 
to the algebraic context (based on the results of Sabbah \cite{Sab}).\\

So we see that the proof of theorem \ref{thm:main} is very simple, once we use the Verdier Riemann-Roch 
theorem (\ref{smoothVerdier}) for smooth morphisms and the specialization property 
(\ref{spec}) of the Chern-class transformation $c_{*}$.
The main point of our paper is therefore just to explain the definition of the specialization functor
$Sp_{X\backslash Y}$ and the "generalized vanishing cycles" $\Phi_{i}$ for constructible functions,
and to show how to use them for the applications.\\

Our "defect-formula"  (\ref{defect})  becomes very simple in the following two cases:

\begin{cor} \label{cor:nonchar}
Consider a regular embedding  $i: X\hookrightarrow Y$ as before and suppose that $\alpha\in F(Y)$
has the property $\Phi_{i}(\alpha) = 0$ (i.e. $\pi^{*}i^{*}(\alpha) = Sp_{X\backslash Y}(\alpha)$).  
Then one has a corresponding Verdier-type Riemann-Roch-theorem for $\alpha$:
\[c_{*}\bigl(i^{*}(\alpha)\bigr) = c^{*}(N_{X}Y)^{-1} \cap i^{*}\bigl(c_{*}(\alpha)\bigr) \;. \]
\end{cor}

The assumption on $\alpha$ in the above corollary is for example satisfied, if $X$ and $Y$ are smooth and 
$X$ is transverse to a complex analytic (or algebraic) Whitney b-regular stratification adapted to $\alpha$
(i.e. $\alpha$ is (locally) constant on the strata of this stratification). But we will not give a proof of this fact
(which uses Fourier-transformation of constructible functions (or sheaves)). Let us just remark that in this 
case we have
\[ c^{*}(N_{X}Y)^{-1} = c^{*}(TX)\cup c^{*}(i^{*}TY)^{-1}\;,\]
and a micro-local proof of the corresponding Verdier-type Riemann-Roch theorem was given in \cite{Sch2}
(based on the theory of (non-characteristic pullback of)Lagrangian cycles).\\

The other case is given by the

\begin{cor} \label{cor:main}
Consider a regular embedding  $i: X\hookrightarrow Y$ of codimension one. 
Then one has for the "defect" in the
corresponding Verdier-type Riemann-Roch-formula the following simple description:
\begin{equation} \label{defect2}
c^{*}(N_{X}Y)^{-1} \cap i^{*}\bigl(c_{*}(\cdot)\bigr)
- c_{*}\bigl(i^{*}(\cdot)\bigr) 
= c^{*}(N_{X}Y)^{-1} \cap c_{*}\bigl(\mu(\cdot)\bigr)\;,
\end{equation} 
with $\mu: F(Y)\to F(X)$ defined as in (\ref{defmu}) .
\end{cor}

Note that this follows from (\ref{defect}) and the definition (\ref{defmu}) of $\mu$
by the following easy calculations:
\[c_{*}\bigl(\Phi_{i}(\cdot)\bigr) = c_{*}\Bigl(\pi^{*}\bigl(\mu(\cdot)\bigr)\Bigr)
- c_{*}\Bigl(k_{*}\bigl(\mu(\cdot)\bigr)\Bigr)\]
\[= c^{*}(T_{\pi}) \cap \pi^{*}\Bigl(c_{*}\bigl(\mu(\cdot)\bigr)\Bigr)
- k_{*}\Bigl(c_{*}\bigl(\mu(\cdot)\bigr)\Bigr) \;.\]
For the last equality we use the Verdier-type Riemann-Roch theorem for the smooth projection
$\pi$ (for the first term), and the functoriality of $c_{*}$ with respect to the proper map
$k$ (for the second term). Then the claim follows from
\[k^{*}\Bigl(c^{*}(T_{\pi}) \cap \pi^{*}\Bigr) - k^{*}k_{*} =
c^{*}(k^{*}T_{\pi}) \cap k^{*}\pi^{*}  - k^{*}k_{*} \]
\[= c^{*}(N_{X}Y) \cap -\; c^{1}(N_{X}Y) \cap = \Bigl(1 + c^{1}(N_{X}Y) -c^{1}(N_{X}Y)\Bigr)\cap
= 1 \cap \;,\]
where we used the self-intersection formula $ k^{*}k_{*} = c^{1}(N_{X}Y)\cap$ (compare \cite[prop.2.6(c), cor.6.3]{Ful} for the algebraic case).\\

If we apply this corollary to $1_{Y}$ for a smooth $Y$, we get as a very special case the formula
\begin{equation} \label{Milnor2}
M_{0}(X) = c^{*}(N_{X}Y)^{-1} \cap c_{*}\bigl(\mu(1_{Y})\bigr)\;.
\end{equation} 
As we will explain later on, the constructible function $\mu(1_{Y})$ is exactly the corresponding
function $\mu$ used in \cite{PP} (up to the sign $(-1)^{dim(Y)-1}$, by which our definition of
$M_{0}(X)$ differs from the convention used in \cite{PP}). So in this special case we get back
the main result of Parusi\'{n}ski-Pragacz \cite[thm.0.2]{PP} giving a nice formula of the
Milnor class of a ("global") hypersurface inside a (compact) complex manifold (which was conjectured
by Yokura). Even in this special case, we think that our more general approach is simpler
than the proof given in \cite{PP}, which is based on deep results about suitable
characteristic cycles. Parusi\'{n}ski-Pragacz used already in \cite{PP} a specialization
argument (i.e. Verdier's specialization result for $c_{*}$)
for the proof of another formula, which for projective $Y$ implies a slightly weaker
version of the above formula for $M_{0}(X)$
(by some lengthy calculation, compare \cite[p.363-367]{Y}).\\

But up to now, this specialization property of $c_{*}$ was not applied to the
specialization map 
\[Sp_{X\backslash Y}: F(Y) \to F_{mon}(C_{X}Y)\;,\]
which is so basic for our approach (and even so natural, if one looks at the definition
of the Gysin map  in homology for regular embeddings). In the special case of the constructible function
$1_{Y}$, we get by theorem \ref{thm:main} the following generalization of
the formula of Parusi\'{n}ski-Pragacz  for a general  regular embedding $i: X\hookrightarrow Y$
of complex "algebraic"  (or compact analytic) spaces, with $Y$ smooth:
\begin{equation} \label{Milnor3}
M_{0}(X) = c^{*}(N_{X}Y)^{-1} \cap k^{*}\Bigl( c_{*}\bigl(\Phi_{i}(1_{Y})\bigr)\Bigr) \;,
\end{equation} 
with $k: X\to N_{X}Y$ the zero-section. The difference to the codimension one case comes from the fact
that in  general  the function $ \Phi_{i}(1_{Y }) \in F_{mon}(N_{X}Y)$ is much more complicated
(compared to the codimension one case)!
Newertheless, in the end of this paper we will explain, that  one can deduce many (known) properties of
$M_{0}(X)$ from this formula (\ref{Milnor3}).
\\

Let us end this introduction with a general "defect-formula" for a local complete intersection 
morphism $f: X\to Y$, which follows directly from the factorization (\ref{factVerdier}),
the Verdier-formula (\ref{smoothVerdier}) for smooth morphisms and theorem \ref{thm:main}:

\begin{cor} \label{cor:main2}
Consider a local complete intersection morphism  $f=p\circ i: X\hookrightarrow Y$ 
with a global factorization into a regular embedding $i: X\to Z$ followed by a smooth morphism
$p: Z\to Y$
of complex "algebraic" (or compact analytic) spaces. Then one has for the "defect" in the
corresponding Verdier-type Riemann-Roch-formula the following explicit description:
\begin{equation} \label{defect3}
c^{*}(T_{f}) \cap f^{*}\bigl(c_{*}(\cdot)\bigr)
- c_{*}\bigl(f^{*}(\cdot)\bigr) 
= c^{*}(N_{X}Z)^{-1} \cap k^{*}\Bigl(c_{*}\bigl(\;\Phi_{i}(p^{*}(\cdot))\;\bigr)\Bigr)\;,
\end{equation}
with 
$c^{*}(T_{f})\cap := \Bigl( c^{*}(i^{*}T_{p})\cup \bigl(c^{*}(N_{X}Z)\bigr)^{-1}\Bigr) 
\cap $.
Here $k: X\hookrightarrow N_{X}Z$ is the "zero-section" of the normal-bundle
$\pi: N_{X}Z\to X$ and $k^{*}$ on homology is the inverse 
of the isomorphism $\pi^{*}: H_{*}(X)\to H_{*}(N_{X}Z)$.
\end{cor}

\section{The Verdier specialization functor for constructible functions}

In this section we explain some basic constructions for constructible functions. Here we work with a complex analytic
space $X$, and a function $\alpha: X\to \bb{Z}$ is called (analytically) constructible, if it satisfies one of
the following two equivalent properties:
\begin{enumerate}
\item $\alpha$ is a locally finite sum $\alpha = \sum_{j} n_{j}\cdot 1_{Z_{j}}$, with $n_{j}\in \bb{Z}$ and
$1_{Z_{j}}$ the indicator function of the closed analytic subset $Z_{j}$ of $X$.
\item   
$\alpha$ is (locally) constant on the strata of a
complex analytic Whitney b-regular stratification of $X$
(which we then call a stratification adapted to $\alpha$).
\end{enumerate}
Note that on a compact analytic space a locally finite sum as above is indeed a finite sum.\\

This notation is closely related to the much more sophisticated language of (analytically) constructible (complexes 
of) sheaves on $X$. A sheaf ${\cal F}$ of vector-spaces on $X$ with finite dimensional stalks
is (analytically) constructible, if there
exists an analytic Whitney b-regular stratification as above such that the restriction of ${\cal F}$ to all
strata is locally constant. Similarly, a bounded complex of sheaves is constructible, if all it cohomology sheaves
have this property, and we denote by $D^{b}_{c}(X)$ the corresponding derived category of bounded constructible
complexes on $X$. Then the Grothendieck group of the triangulated category $D^{b}_{c}(X)$ is isomorphic
to the Grothendieck group of the abelian category of constructible sheaves (by associating to a constructible
complex the alternating sum of its cohomology sheaves), and we denote it by $K_{const}(X)$.\\

Since we assume that all stalks of a constructible complex are finite dimensional, we get by taking stalkwise
the Euler characteristic a natural group 
homomorphism $\chi_{X}: K_{const}(X)\to F(X)$ into the group $F(X)$ of constructible functions (on $X$).
And it is very easy to show that $\chi_{X}$ is surjective.\\

If $X$ is the analytic space associated to a seperated scheme of finite type over $spec(\bb{C})$, then one has
the corresponding notation of an algebraically constructible (complex of) sheaves or function on $X$
(by using an algebraic stratification of $X$), and in this case any algebraically constructible function is just
a finite linear combination of functions $1_{Z}$ with $Z$ a closed algebraic subset of $X$.
As before one gets a corresponding group epimorphism $\chi_{X}$, and in the following we treat the analytic and algebraic
context (as far as possible) at the same time.\\

As is well known (and explained in detail in \cite{Sch}), all the usual functors in sheaf theory, which respect
the corresponding category of constructible complexes of sheaves, induce by the epimorphism $\chi_{X}$ well-defined
group homomorphisms on the level of constructible functions. We just recall those, which are important for later
applications or definitions. Let $f: X \to Y$ be a holomorphic (or algebraic) map of complex spaces (associated
to seperated schemes of finite type over $spec(\bb{C})$). Then one has the following transformations:
\begin{itemize}
\item {\bf pullback:} $f^{*}: F(Y)\to F(X);\; \alpha\mapsto \alpha\circ f$, which corresponds to the usual pullback
of sheaves.
\item {\bf Euler characteristic:} Suppose $X$ is compact and $Y=\{pt\}$ is a point. Then one has 
$\chi: F(X)\to \bb{Z}$,
corresponding to $R\Gamma(X,\cdot)$ on the level of constructible complexes of sheaves. By linearity it is  
characterized by the convention, that $\chi(1_{Z}) $ for a closed analytic (algebraic) subspace $Z$ of $X$ 
is just the usual Euler characteristic of the (finite dimensional) cohomology $H^{*}(Z)$ of $Z$.
\item {\bf proper pushdown:} Suppose $f$ is proper. Then one has $f_{*}: F(X)\to F(Y)$, corresponding to $Rf_{*}$
on the level of constructible complexes of sheaves. Explicitly it is given by $f_{*}(\alpha)(y) := 
\chi(\alpha |\{f=y\})$, and in this form it goes back to the paper \cite{M} of MacPherson
(Note that $\{f=y\}$ is compact, since $f$ is proper so that we can use the Euler characteristic $\chi$).
\item {\bf nearby cycles:} Assume $Y=\bb{C}$ and let $X_{0}:=\{f=0\}$ be the zero fiber. Then one has
$\psi_{f}: F(X)\to F(X_{0})$, corresponding to Deligne's nearby cycle functor.
This was first introduced in \cite[Exp.7,prop.3.4]{CEP} by using resolution of singularities
(compare with \cite{Sch} for another approach using stratification theory).
By linearity, $\psi_{f}$ is uniquely defined by the convention that for a closed analytic (algebraic)
subspace $\psi_{f}(1_{Z})(x)$ is just the Euler-characteristic of a local Milnor fiber $M_{f|Z,x}$ of $f|Z$ at $x$
\cite[Exp.7,prop.4.1]{CEP}. Here this local Milnor fiber at $x$ is given  
by $M_{f|Z,x}:=Z\cap B_{\epsilon}(x)\cap \{f=y\}$, with $0<y<<\epsilon<<1$ and $B_{\epsilon}(x)$ an open
(or closed) ball of radius $\epsilon$ around $x$ (in some local coordinates). 
\item {\bf vanishing cycles:} Assume $Y=\bb{C}$ and let $i: X_{0}:=\{f=0\}\hookrightarrow X$ be the inclusion
of the zero-fiber. Then one has
$\phi_{f}: F(X)\to F(X_{0});\;\phi_{f} := \psi_{f} - i^{*} $, corresponding to Deligne's vanishing cycle functor.
\end{itemize}

Using these notations, we can now introduce the "constructible function version" of Verdier's specialization 
functor. Let $i: X\hookrightarrow Y$ be a closed embedding of complex analytic spaces
(associated to seperated schemes of finite type over $spec(\bb{C})$). Then the deformation to the normal cone
$C_{X}Y$ of $X$ in $Y$ is a complex analytic ("algebraic") space $M$ together with two morphisms
$\pi: M\to Y$ and $h: M\to \bb{C}$, with $h$ flat so that one has a cartesian diagram
\begin{equation} \label{eq:defnormal} \begin{CD}
C_{X}Y @>s >> M @< j << Y \times \bb{C}^{*}\\
@VVV  @VV h V  @VVpr V\\
\{0\} @>>> \bb{C} @<<< \bb{C}^{*}\;.
\end{CD} \end{equation}

For details of its construction we refer to \cite[Ch.5]{Ful} and \cite[Exp.IX]{GA}.
Let us only recall that $M=\tilde{M}\backslash \tilde{Y}$, with $\tilde{M}$ the blow-up of $Y\times \bb{C}$
along $X\times \{0\}$ and $\tilde{Y}$ the blow-up of $Y$ along $X$ (suitable embedded into $\tilde{M}$)
so that $\pi$ and $h$ are induced from the corresponding projections by the blow-up map $b$:
\begin{equation} \label{eq:b} \begin{CD}
M = \tilde{M}\backslash \tilde{Y} \hookrightarrow \tilde{M} @> b >> Y\times \bb{C} \;.
\end{CD} \end{equation}

Then $
Sp_{X\backslash Y}: F(Y)\to F_{mon}(C_{X}Y)$ is defined as $Sp_{X\backslash Y}:= \psi_{h}\circ \pi^{*}$.
This corresponds exactly to Verdier's definition \cite[p.352,p.358]{V} for constructible sheaves,
since $\psi_{h}(p^{*}{\cal F})$ depends by definition only on 
$j^{*}\pi^{*}{\cal F} = (\pi\circ j)^{*}{\cal F}$, and $\pi\circ j$ is by construction
of the deformation just the projection $ Y \times \bb{C}^{*} \to Y$. Note that the construction of $(M,\pi,h)$ is
functorial in the pair $(Y,X)$.\\

By \cite{V} we therefore get the following important properties of the functor $Sp_{X\backslash Y}$
(where the property (SP4) related to duality is trivial for constructible functions, so that we omit it):
\begin{itemize}
\item[(SP0)] {\bf localization:} $Sp_{X\backslash Y}$ commutes with restriction to open subsets. 
\item[(SP1)]  {\bf monodromy:} $Sp_{X\backslash Y}$ maps to $F_{mon}(C_{X}Y)$.
\item[(SP2)] {\bf proper direct image:} Consider a cartesian diagram
\[ \begin{CD} 
X' @>i'>> Y'\\
@VVV  @VVf V\\
X @>> i> Y
\end{CD} \]
with $f$ proper. Then also the induced map $C(f): C_{X'}Y'\to C_{X}Y$ is proper with
\begin{equation} 
 Sp_{X\backslash Y}\circ f_{*} = C(f)_{*}\circ Sp_{X'\backslash Y'}\;. 
\end{equation}
\item[(SP3)] {\bf smooth base change:}  Consider a cartesian diagram
\[ \begin{CD} 
X' @>i'>> Y'\\
@VVV  @VVf V\\
X @>> i> Y
\end{CD} \]
with $f$ smooth. Then also the induced map $C(f): C_{X'}Y'\to C_{X}Y$ is smooth with
\begin{equation} 
C(f)^{*}\circ Sp_{X\backslash Y} =  Sp_{X'\backslash Y'}\circ f^{*} \;. 
\end{equation}
\item[(SP5)] {\bf restriction to the vertex:} $Sp_{X\backslash Y}(\alpha)|X = \alpha|X$ for all $\alpha\in F(Y)$.
\item[(SP6)] {\bf normalisation:} Suppose $X$ is defined on $Y$ by one equation $\{g=0\}$.
Then $\pi\times C(g): C_{X}Y = N_{X}Y\to X\times \bb{C}$ is an isomorphism. 
Let $s: X\to  N_{X}Y$ be the section corresponding
to $id\times \{1\}$ under this isomorphism. Then $s^{*}\circ Sp_{X\backslash Y} = \psi_{g}$.
\end{itemize}

Recall the we introduced in (\ref{vanish}) the "generalized vanishing cycles"
$\Phi_{i}$ of $i$ as the transformation
\[\Phi_{i}:= Sp_{X\backslash Y}
- \pi^{*}i^{*} 
: F(Y) \to F_{mon}(C_{X}Y)\;.\]
For a cartesian diagram as in (SP2,3) one gets also the commutative diagram (with $k,k'$ the "zero-sections"):
\begin{equation} \label{flatsp} \begin{CD} 
X' @>k'>> C_{X'}Y' @>\pi' >> X'\\
@Vf VV  @VV C(f) V  @VV f V\\
X @>> k> C_{X}Y @>> \pi > X\;,
\end{CD} \end{equation}
which is cartesian for a flat morphism $f$.\\

Therefore one gets similar properties for $\Phi_{i}$:
\begin{itemize}
\item[(gV0)] {\bf localization:} $\Phi_{i}$ commutes with restriction to open subsets. 
\item[(gV1)]  {\bf monodromy:} $\Phi_{i}$ maps to $F_{mon}(C_{X}Y)$.
\item[(gV2)] {\bf proper direct image:} Consider a cartesian diagram
\[ \begin{CD} 
X' @>i'>> Y'\\
@VVV  @VVf V\\
X @>> i> Y
\end{CD} \]
with $f$ flat and proper
(or more generally with $f$ proper such that (\ref{flatsp}) is cartesian). 
Then also the induced map $C(f): C_{X'}Y'\to C_{X}Y$ is proper with
\begin{equation} 
 \Phi_{i}\circ f_{*} = C(f)_{*}\circ \Phi_{i'}\;. 
\end{equation}
\item[(gV3)] {\bf smooth base change:}  Consider a cartesian diagram
\[ \begin{CD} 
X' @>i'>> Y'\\
@VVV  @VVf V\\
X @>> i> Y
\end{CD} \]
with $f$ smooth. Then also the induced map $C(f): C_{X'}Y'\to C_{X}Y$ is smooth with
\begin{equation} 
C(f)^{*}\circ \Phi_{i} = \Phi_{i'} \circ f^{*} \;. 
\end{equation}
\item[(gV5)] {\bf restriction to the vertex:} $\Phi_{i}(\alpha)|X = 0$ for all $\alpha\in F(Y)$.
\item[(gV6)] {\bf normalisation:} Suppose $X$ is defined on $Y$ by one equation $\{g=0\}$.
Then $\pi\times C(g): C_{X}Y = N_{X}Y\to X\times \bb{C}$ is an isomorphism. 
Let $s: X\to  N_{X}Y$ be the section corresponding
to $id\times \{1\}$ under this isomorphism. Then $s^{*}\circ \Phi_{i} = \phi_{g}$.
\end{itemize}

Let us now explain the definition of the corresponding specialization map $sp$ in homology.
First we consider the analytic context (compare \cite[p.217,218]{GA} and \cite[p.149]{CEP}).
Since $\pi: M\backslash C_{X}Y = Y\times \bb{C}^{*} \to \bb{C}^{*}$ is just the projection,
one gets that $(R\pi_{!}\bb{Z}_{M})|\bb{C}^{*}$ is isomorphic to the constant complex $R_{c}(Y,\bb{Z})$.
Therefore one gets a generalization morphism (independent of the choice of $t\in \bb{C}^{*} $):
\[ R\Gamma_{c}(C_{X}Y,\bb{Z})\simeq (R\pi_{!}\bb{Z}_{M})_{0} \simeq R(\bb{C},R\pi_{!}\bb{Z}_{M}) 
\to (R\pi_{!}\bb{Z}_{M})_{t}\simeq R\Gamma_{c}(Y,\bb{Z})\;,\]
which induces by duality a specialization homomorphism in Borel-Moore homology:
\begin{equation} \label{specialhom}  \begin{CD}
H^{BM}_{*}(C_{X}Y,\bb{Z})  @<sp<< H^{BM}_{*}(Y,\bb{Z})\;.
\end{CD} \end{equation}
The definition in the algebraic context is simpler (compare \cite[ch.5]{Ful} and \cite[p.198]{GA}):
\begin{equation} \label{specialchow}  
sp:=s^{*}\circ \pi^{*}:   A_{*}(Y,\bb{Z})\to  A_{*}(C_{X}Y,\bb{Z}) \;.
 \end{equation}
Here $\pi^{*}$ is the flat pullback and $s^{*}$ is the corresponding Gysin transformation for divisors.
The algebraic and topological definitions commute with the cycle map $cl: A_{*}(\cdot)\to  H^{BM}_{*}(\cdot,\bb{Z})$
(compare \cite[p.222,cor.8.10]{GA}).\\

Now we can formulate the important specialization property for the Chern-Schwartz-MacPherson classes 
$c_{*}$ (Compare with \cite[p.352]{Ful} and \cite[thm.6.5,p.211]{GA}
for the corresponding specialization property of
the "Todd-homology transformation" $\tau_{*}$ of Baum-Fulton-MacPherson).

\begin{thm} \label{spezialc*}
Consider the complex algebraic context with algebraic constructible functions $F(\cdot)$
 and $H_{*}(\cdot)$ one of 
the two possible homology theories
(i.e. the Chow-homology $A_{*}(\cdot)$ or the Borel-Moore homology $H^{BM}_{*}(\cdot,\bb{Z})$),
or the complex analytic context with $Y$ compact. Then the following diagram is commutative:
\begin{equation} \begin{CD}
F(Y) @> Sp_{X\backslash Y} >> F_{mon}(C_{X}Y)\\
@V c_{*} VV   @VV c_{*} V \\
H_{*}(Y) @>> sp > H_{*}(C_{X}Y) \;.
\end{CD} \end{equation}
\end{thm}

\begin{proof} The analytic and algebraic case with Borel-Moore homology follows directly from
a corresponding result of Verdier for the nearby cycle functor $\psi_{h}$ \cite[p.157,thm.5.1]{CEP},
with $h$ the map of the deformation to the normal cone (\ref{eq:defnormal}):
\[c_{*}\Bigl(\psi_{h}(\pi^{*}(\alpha))\Bigr) = sp\Bigl(c_{*}((\pi^{*}\alpha)|Y\times \{t\})\Bigr)\]
(for $0<t<<1$). Of course $c_{*}((\pi^{*}\alpha)|Y\times \{t\})$ corresponds to $c_{*}(\alpha)$ under the
identification $ Y\times \{t\}\simeq Y$.
Note that he assumes that the "situation is compactifiable", which in the algebraic context is automatically
the case. In the analytic context it follows from the compactness of $Y$ and the construction of the deformation
to the normal cone (e.g. $\pi^{*}(\alpha) = (\tilde{\pi}^{*}(\alpha))|M$, with $\tilde{\pi}:
\tilde{M}\to Y$ the induced map, and the extension $\tilde{M}\to Y\times \bb{C}\to \bb{C}$ of $h$
is proper for $Y$ compact).

In the algebraic context with the Chow-group as our homology theory, the claim follows from
a similar result of Kennedy for the functor $\psi_{h}$ \cite[thm.2)]{Ken2} (or from the fact that
the argument of Verdier also works on the level of Chow-groups).
Note that Kennedy uses the Lagrangian approach so that he works with spaces that can be embedded
into manifolds. Moreover one has to use the result of Sabbah \cite{Sab} (compare also with 
\cite{BMM,Fu,Gi,Sch3}), 
that the used specialization of Lagrangian cycles corresponds over the ground field $\bb{C}$ to the nearby 
cycle functor for constructible functions. For a general constructible function one can use resolution
of singularities (and induction on the dimension of the support of the constructible function),
to reduce to the case of a function on a smooth space. This works because all considered transformations
commute with proper push down:
\begin{itemize}
\item For the Chern-Schwartz-MacPherson class $c_{*}$ this is one of its defining properties.
\item For the nearby cycle functor of constructible functions this follows from a similar statement
for constructible sheaves (and this is exactly what is used for property (SP2) of the specialization
 functor): $p_{*}\circ \psi_{h\circ p} = \psi_{h}\circ p_{*}$ for a proper morphism $p$. 
\end{itemize}
\end{proof}

\begin{rem}
We believe that theorem \ref{spezialc*} is true in the general analytic context, without any compactness
assumptions, but for this one has to find a new proof. If this can be done, then all results of this paper extend 
without any modification to this more general analytic context (without any compactness assumptions).
By the {\bf localization} property (SP0) this is for example true, if one restricts the
complex algebraic (or compact analytic) context to an open subset.
\end{rem}

Theorem \ref{spezialc*} implies directly our important formula (\ref{sp}) and therefore (by the other
arguments of the introduction) also our main theorem \ref{thm:main}. This follows from the definition of
the Gysin homomorphism $i^{*}: H_{*}(Y)\to H_{*}(X)$ for a regular embedding $i: X\hookrightarrow Y$, which by
definition (\cite[Ch.5,Ch.6]{Ful} or \cite[Exp.IX]{GA}) is just the composition $i^{*}:= k^{*}\circ sp$.
Here $k: X\hookrightarrow C_{X}Y=N_{X}Y$ is the "zero-section" of the normal-bundle
$\pi: N_{X}Y\to X$ (which is a vector-bundle for a regular embeddding) 
and
$k^{*}$ on homology is the inverse 
of the isomorphism $\pi^{*}: H_{*}(X)\to H_{*}(N_{X}Y)$.

\section{Milnor-classes}

In this last section, we explain some simple facts which follow from our generalization (\ref{Milnor3}) of
the formula of Parusi\'{n}ski-Pragacz  for a general  regular embedding $i: X\hookrightarrow Y$
of complex algebraic  (or compact analytic) spaces, with $Y$ smooth:
 
\[M_{0}(X) = c^{*}(N_{X}Y)^{-1} \cap k^{*}\Bigl( c_{*}\bigl(\Phi_{i}(1_{Y})\bigr)\Bigr) \;,
\] 
with $k: X\to N_{X}Y$ the zero-section. \\

We already explained in the introduction how to get in the codimension one case the
main result of Parusi\'{n}ski-Pragacz as formulated in \cite[thm.0.2]{PP}.
Moreover, the definition of the {\bf nearby cycles} for constructible functions in terms
of local Milnor fibers and the normalization properties (SP6) and (gV6) show, that our
constructible function $(-1)^{dim(Y)-1}\cdot\mu(1_{Y})$ is exactly the function $\mu$ used in \cite{PP}.
Note that most authors include a similar sign-factor into their definition of the Milnor-class.
But we prefer our definition because of its relation to the "(generalized) vanishing cycles" of
constructible functions.\\

The next important fact is 
\begin{equation}
supp(\Phi_{i}(1_{Y}))\subset N_{X}Y|X_{sing} \;,
\end{equation} 
with $X_{sing}$ the singular
locus of $X$. This follows from the {\bf localization} property (gV0) of $\Phi_{i}$, and the fact 
that for a regular embedding $X\hookrightarrow Y$ of a smooth subspace the deformation to the normal
cone $M$ and the induced map $h: M\to \bb{C}$ are also smooth! This can be easily checked in local 
coordinates, and is of course very important for the D-module approach to the specialization
functor (compare \cite{La}). But for a smooth morphism $h$ one has $\psi_{h}(1_{M}) = 1_{N_{X}Y}$
and therefore $\Phi_{i}(1_{M})=0$.\\

Let us denote by $i_{\alpha}: S_{\alpha}\hookrightarrow X$ the (closed) inclusions of the
connected components of $X_{sing}$, and by $n_{\alpha}: N_{X}Y|S_{\alpha} \to N_{X}Y$
the corresponding inclusions.
Then we get $\Phi_{i}(1_{Y}) = \sum_{\alpha} \;n_{\alpha *} (\mu_{\alpha})$, with
$\mu_{\alpha} := n_{\alpha}^{*} \Phi_{i}(1_{Y})$ and
\[c_{*}(\Phi_{i}(1_{Y})) = \sum_{\alpha} \;n_{\alpha *} (c_{*}(\mu_{\alpha}))\;.\]
>From the cartesian diagram (with $k$ the zero-section)
\[ \begin{CD}
S_{\alpha} @> s_{\alpha}  >> X \\
@V k_{\alpha} VV  @VV k V \\
N_{X}Y|S_{\alpha} @>> n_{\alpha} > N_{X}Y 
\end{CD} \]
one gets by the base-change property $k^{*}\circ n_{\alpha *} = s_{\alpha *}\circ k_{\alpha}^{*}$ 
the following localization formula for the Milnor-class:
\begin{equation} 
M_{0}(X) = \sum_{\alpha}\; c^{*}(N_{X}Y)^{-1} \cap 
s_{\alpha *}k_{\alpha}^{*}\Bigl(c_{*}(\mu_{\alpha})\Bigr)\]
\[= \sum_{\alpha}\; s_{\alpha *} \Bigl(\; c^{*}(N_{X}Y|S_{\alpha})^{-1} \cap 
 k_{\alpha}^{*}(c_{*}(\mu_{\alpha}))\;\Bigr)\;.
\end{equation}

So we get a localization result similar to \cite[thm. A]{BLSS} or \cite[thm.5.2]{BLSS2},
and it would be very interesting to compare our "localized Milnor-class"
\begin{equation} \label{eq:locMil}
M_{0}(X,S_{\alpha}) :=   c^{*}(N_{X}Y|S_{\alpha})^{-1} \cap 
k_{\alpha}^{*}(\;c_{*}(\mu_{\alpha})\;) \in \; H_{*}(S_{\alpha})
\end{equation}
with the corresponding notation introduced in \cite{BLSS,BLSS2} with the help of obstruction theory.\\

Of course these "localized Milnor-classes" $M_{0}(X,S_{\alpha})$ are very difficult to calculate in general.
But by dimension reasons, the corresponding multiplicity in $\bb{Z} \cdot [S'_{\alpha}] \subset 
H_{*}(S_{\alpha})$, with $[S'_{\alpha}]$ the fundamental-class of an irreducible component $S'_{\alpha}$
of $S_{\alpha}$ is nothing but the "generic value" of $\mu_{\alpha}$ (or $\Phi_{i}(1_{Y})$)
 on $N_{X}Y|S'_{\alpha}$.
Note that the restriction of the constructible function $\mu_{\alpha}$ to the irreducible
space $N':=N_{X}Y|S'_{\alpha}$ has a "generic value" (i.e. the coefficient of $1_{N'}$ in a description
of this constructible function as a linear combination of indicator-functions of closed irreducible
subsets of $N'$). \\

Especially if $S_{\alpha}=\{x_{\alpha}\}$ is just a point, 
then it follows from (the corresponding local Milnor-fibration and)
the definition of $\Phi_{i}(1_{Y})$ in terms of nearby-cycles that
this "generic value" is just $\chi(M_{x_{\alpha}})-1$, with $M_{x_{\alpha}}$ a local
Milnor-fiber of the local complete intersection $X$ at the isolated singular point $x_{\alpha}$.
But this Milnor-fiber is wedge of spheres of dimension $d:=dim_{x_{\alpha}}(X)$ so that
$(-1)^{d}\cdot(\chi(M_{x_{\alpha}})-1)$ is just the usual Milnor-number of this isolated
local complete intersection singularity (i.e. the corresponding number of spheres). 
So for a complete intersection $X$ with only isolated singularities we get back a result of Suwa 
\cite[thm.,p.68]{Su}
that the Milnor-class $M_{0}(X)$ is just the sum
\begin{equation}
M_{0}(X) = \sum \; (\chi(M_{x_{\alpha}})-1)\cdot s_{\alpha *}([\{x_{\alpha}\}]) 
\end{equation}
of the classes of the singular points $x_{\alpha}$ weighted by the corresponding
Milnor-numbers (up to suitable signs). \\

There is also another special case, where one can give a simpler formula for the 
"localized Milnor-class" $M_{0}(X,S_{\alpha})$. 
Suppose $\mu_{\alpha} = \sum \;m_{j}\cdot 1_{V_{j}}$, with $p_{j}: V_{j}\to S_{j}$ a closed 
subvector-bundle of $N_{X}Y|S_{j}$ for a closed subspace $S_{j}\subset S_{\alpha}$. 
By linearity (and to simplify the notations) we can assume $\mu_{\alpha}=1_{V}$
with $p: V\to S$ a closed subvector-bundle of $N_{X}Y|S$ for a closed subspace $S\subset S_{\alpha}$.
Consider the commutative diagram
\[ \begin{CD}
V @> j >>  N_{X}Y|S  @> i'   >> N_{X}Y|S_{\alpha} \\
@A \kappa AA  @AA \kappa' A  @AA k_{\alpha} A \\
S @= S @>> i > S_{\alpha} \;,
\end{CD} \]
where the vertical maps are the "zero-sections". Then one gets:
\[k_{\alpha}^{*}(c_{*}(1_{V})) = k_{\alpha}^{*}i'_{*}(c_{*}(1_{V})) = i_{*}\kappa'^{*}(c_{*}(1_{V}))\]
\[= i_{*}\kappa^{*}j^{*}j_{*}(c_{*}(1_{V})) = i_{*}\kappa^{*}(c^{d}(N)\cap c_{*}(1_{V}))
= i_{*}\Bigl(c^{d}(\kappa^{*}N)\cap \kappa^{*}(c_{*}(1_{V}))\Bigr)\;.\]
Here we consider the constructible function $1_{V}$ and its chern-class $c_{*}(1_{V})$  
on the corresponding spaces. Moreover, we use for the inclusion $j$ the "self-intersection formula"
with $N$ the corresponding normal-bundle of rang $d$ so that 
$\kappa^{*}N$ is isomorphic to the quotient-bundle $(N_{X}Y|S)/V$.
Since $\kappa^{*}(c_{*}(1_{V})) = c^{*}(V)\cap c_{*}(1_{S})$ (by the Verdier-formula for the smooth
projection $p$), we get altogether
\[k_{\alpha}^{*}(c_{*}(1_{V})) =
 i_{*}\Bigl(\;c^{d}((N_{X}Y|S)/V)\cap c^{*}(V)\cap c_{*}(1_{S})\;\Bigr)\;.\]
By (\ref{eq:locMil}) we get a similar formula for the "localized Milnor-class",
and arguments of this type can be used to relate our general formula (\ref{Milnor3}) to
the product-formula for Milnor-classes obtained in \cite{OY}.\\

Let us end this paper with some functoriality results for the "Milnor-class",
which are easy applications of 
our formula (\ref{Milnor3}), the general results about our "generalized vanishing cycles" $\Phi_{i}$
and known properties of the Chern-Schwartz-MacPherson classes $c_{*}$.\\

First we explain the behaviour of $M_{0}(X)$ under smooth pullback (compare \cite[thm.2.2]{Y3}).
Consider a cartesian diagram
\[ \begin{CD}
X' @> i' >>  Y' \\
@V f VV  @VV f' V \\
X @>> i > Y \;,
\end{CD} \]
with $f'$ smooth and $i: X\hookrightarrow Y$ a closed embedding.
Then the corresponding diagram (\ref{flatsp}) is cartesian, with $f,C(f)$ smooth and
$i': X'\hookrightarrow Y'$ is also a closed embedding with $f^{*}(C_{X}Y) \simeq C_{X'}Y'$. 
Then we get by (gV3) and the Verdier-formula (\ref{smoothVerdier}) for the smooth map $C(f)$:
\[k'^{*}\Bigl(\;c_{*}(\Phi_{i'}(1_{Y'}))\;\Bigr) = 
k'^{*}\Bigl(\;c_{*}(\;(\Phi_{i'}\circ f^{*})(1_{Y})\;)\;\Bigr) =\]
\[k'^{*}\Bigl(\;c^{*}(T_{C(f)})\cap \;C(f)^{*}(\;c_{*}(\Phi_{i}(1_{Y}))\;)\;\Bigr)
= c^{*}(T_{f})\cap \;f^{*}\Bigl(\;k^{*}(c_{*}(\Phi_{i}(1_{Y})))\;\Bigr)\;.\]

Assume we consider the algebraic (or analytic) context with $Y,Y'$ (compact) manifolds
and $i$ a regular embedding. Then $i'$ is also a regular embedding with $f^{*}(N_{X}Y) \simeq N_{X'}Y'$.
Then we get by this calculation and our formula (\ref{Milnor3}) the following
relation between the corresponding "Milnor-classes":
\begin{equation}
M_{0}(X') = c^{*}(T_{f})\cap \;f^{*}M_{0}(X) \;.
\end{equation}

A similar argument applies also to the "proper pushdown" of "Milnor-classes"
studied in \cite{Y3}.  Consider a cartesian diagram as before, but this time with $f'$ proper
(but not necessarily smooth). Assume the closed embeddings $i,i'$ are
regular of the same (local) codimension. Then one has   $f^{*}(N_{X}Y) \simeq N_{X'}Y'$
(since  $f^{*}(N_{X}Y)$ is a subbundle of $N_{X'}Y'$ of the same rang) so that the diagram
(\ref{flatsp}) is also cartesian in this case, and  $C(f)$ is a proper by "base-change". 
By (gV2) and the functoriality of $c_{*}$ with respect to proper morphisms one gets
\[f_{*}k'^{*}\Bigl(\;c_{*}(\Phi_{i'}(1_{Y'}))\;\Bigr) = 
k^{*}C(f)_{*}\Bigl(\;c_{*}(\Phi_{i'}(1_{Y'}))\;\Bigr) =
k^{*}\Bigl(c_{*}(\;\Phi_{i}(f'_{*}(1_{Y'}))\;)\Bigr)\;.\]
  
Assume we consider the algebraic (or analytic) context with $Y,Y'$ (compact) manifolds
and (for simplicity) $Y$ connected. Let $\chi$ be the "generic value" of the constructible
function $f'_{*}(1_{Y'})$ (i.e. the Euler-characteristic of a "generic" fiber of $f'$) so that
$\alpha:=f'_{*}(1_{Y'}) - \chi\cdot 1_{Y}$ is a constructible function of smaller support.
Then one gets by this calculation, the projection-formula and our formula (\ref{Milnor3}) the following
relation between the corresponding "Milnor-classes":
\begin{equation}
 f_{*}M_{0}(X') = \chi\cdot M_{0}(X) + 
 c^{*}(N_{X}Y)^{-1} \cap k^{*}\Bigl( c_{*}\bigl(\Phi_{i}(\alpha)\bigr)\Bigr) \;.
\end{equation}
Especially, if all fibers of $f'$ have the same Euler-characteristic (e.g. $f'$ is smooth
or more generally an "Euler-morphism" as in \cite[cor.3.8]{Y3}), then one gets
\[f_{*}M_{0}(X') = \chi\cdot M_{0}(X) \;.\]

\end{document}